\documentclass[12pt]{article}

\let\oldmarginpar\marginpar
\renewcommand\marginpar[1]{\-\oldmarginpar[\raggedleft\footnotesize #1]%
{\raggedright\footnotesize #1}}

\usepackage{mathptmx}
\usepackage{amsmath}
\usepackage{amscd}
\usepackage{amssymb}
\usepackage{amsthm}
\usepackage{xspace}
\usepackage[all,tips]{xy}
\usepackage[dvips]{graphicx}
\usepackage{verbatim}
\usepackage{syntonly}
\usepackage{hyperref}

\theoremstyle{plain}
\newtheorem{thm}{Theorem}[section]
\newtheorem{cor}[thm]{Corollary}
\newtheorem{prop}[thm]{Proposition}
\newtheorem{lemma}[thm]{Lemma}

\theoremstyle{definition}

\newenvironment{pf}
{\begin{proof}} {\end{proof}}

 \DeclareMathOperator{\Out}{Out}
\DeclareMathOperator{\Mod}{Mod} \DeclareMathOperator{\Cl}{Cl}
\DeclareMathOperator{\SL}{SL} 
\DeclareMathOperator{\GL}{GL}

\DeclareMathOperator{\FI}{FI}\DeclareMathOperator{\NFI}{NFI}
\DeclareMathOperator{\Ave}{Ave}\DeclareMathOperator{\GCD}{gcd}
\DeclareMathOperator{\lcm}{lcm}\DeclareMathOperator{\Eff}{Eff}

\DeclareMathOperator{\D}{D}

\newcommand{\eps}{\varepsilon}
\newcommand{\vp}{\varphi}

\newcommand{\ol}{\overline}

\newcommand{\nid}{\noindent}
\newcommand{\wh}{\widehat}

\newcommand{\iny}{\infty}

\newcommand{\es}{\emptyset}
\newcommand{\co}{\ensuremath{\colon}}

\newcommand{\innp}[1]{\left< #1 \right>}
\newcommand{\abs}[1]{\left\vert#1\right\vert}
\newcommand{\set}[1]{\left\{#1\right\}}
\newcommand{\brac}[1]{\left[#1\right]}
\newcommand{\pr}[1]{\left( #1 \right) }
\newcommand{\su}{\subset}

\newcommand{\bu}{\bigcup}
\newcommand{\ba}{\bigcap}
\newcommand{\bop}{\bigoplus}

\newcommand{\lra}{\longrightarrow}

\newcommand{\B}[1]{\ensuremath{\mathbf{#1}}}

\newcommand{\Cal}[1]{\ensuremath{\mathcal{#1}}}
\newcommand{\Fr}[1]{\ensuremath{\mathfrak{#1}}}

\newcommand{\N}{\ensuremath{\B{N}}}
\newcommand{\Q}{\ensuremath{\B{Q}}}
\newcommand{\R}{\ensuremath{\B{R}}}
\newcommand{\Z}{\ensuremath{\B{Z}}}
\newcommand{\C}{\ensuremath{\B{C}}}

\usepackage{fancyhdr}

\fancypagestyle{plain}

\begin{document}


\title{\textbf{Bertrand's postulate and subgroup growth}}
\author{K. Bou-Rabee and D. B. McReynolds}
\maketitle


\begin{abstract}
\nid In this article we investigate the $L^1$--norm of certain functions on groups called divisibility functions. Using these functions, their connection to residual finiteness, and integration theory on profinite groups, we define the residual average of a finitely generated group. One of  the main results in this article is the finiteness of residual averages on finitely generated linear groups. Whether or not the residual average is finite depends on growth rates of indices of finite index subgroups. Our results on index growth rates are analogous to results on gaps between primes, and provide a variant of the subgroup growth function, which may be of independent interest.
\end{abstract}

\nid keywords: \emph{Bertrand's postulate, residual finiteness, subgroup growth}\smallskip

\nid MSC code: 20E07, 20E18

\section{Introduction and main results}

\nid The study of gaps between successive primes has been a central topic in number theory for more than a hundred years. One classical result known as \emph{Bertrand's postulate} asserts that for any positive integer $n$, there exists a prime between $n$ and $2n$. This assertion was conjectured by Bertrand in 1845 and proved by Chebyshev in 1850 (see \cite{Tch}). The stronger Legendre conjecture asserts a maximum gap of $\sqrt{n}$ (see the related \cite{I}), while the Prime Number Theorem provides many primes in the interval $[n,2n]$.\smallskip\smallskip

\nid One of the main purposes of this article is to investigate related problems for finitely generated groups. The role of primes in our setting will be played by relatively prime families of finite index subgroups $\set{\Delta_j}$ or, more specifically, the indices of such subgroups. Here, we say a family in $\Gamma$ is \emph{relatively prime} if for any distinct pair $\Delta_j,\Delta_k$, we have $\Gamma=\Delta_j\Delta_k$; we will work with a stronger property called \emph{prime} where a Chinese Remainder Theorem holds (see Section \ref{indexGrowthAndConvergenceSection} for a precise definition). Our interest will be in proving results in the same vein as Bertrand's postulate.\smallskip\smallskip

\nid \textbf{Question 1.} (Bertrand's postulate; classical) \emph{On a finitely generated group $\Gamma$, when does there exist a prime (or normal) family of finite index subgroups $\set{\Delta_j}$ and a constant $d>1$ such that
\[ [\Gamma:\Delta_j] \leq [\Gamma:\Delta_{j+1}] \leq d[\Gamma:\Delta_j]? \] }

\nid In the case $\Gamma = \Z$ and the prime family of subgroups is $\set{p\Z}$, an affirmative answer to Question 1 is given by Bertrand's postulate. We could be less ambitious and allow ourselves small powers.\smallskip\smallskip

\nid \textbf{Question 2.} (Bertrand's postulate; small powers) \emph{On a finitely generated group $\Gamma$, when does there exist a prime (or normal) family of finite index subgroups $\set{\Delta_j}$ and a constant $\delta>0$ such that
\[ [\Gamma:\Delta_j] \leq [\Gamma:\Delta_{j+1}] \leq [\Gamma:\Delta_j]^{1+\delta}? \] }

\nid Our first result resolves Question 1 for finitely generated linear groups.

\begin{thm}\label{BertrandForArithmeticLattices}
Let $\Gamma$ be an infinite finitely generated linear group over $\C$. Then there exists a constant $d$ and an infinite family of finite index subgroups $\set{\Delta_j}$ such that
\[ [\Gamma:\Delta_j] < [\Gamma:\Delta_{j+1}] \leq d[\Gamma:\Delta_j]. \]
In addition, there exists a finite index subgroup $\Gamma_0$ of $\Gamma$ such that $\set{\Delta_j}$ is a normal, prime family in $\Gamma_0$.
\end{thm}

\nid The proof of Theorem \ref{BertrandForArithmeticLattices} uses the Lubotzky Alternative, the Strong Approximation Theorem, and Bertrand's postulate. We also prove the existence of families that resolve Question 2 that avoids the use of the Lubotzky Alternative and the Strong Approximation Theorem.

\begin{thm}\label{PowerBertrand}
Let $\Gamma$ be an infinite finitely generated linear group over $\C$. Then for any $\delta>0$, there exists an infinite normal family of finite index subgroups $\set{\Delta_j}$ such that
\[ [\Gamma:\Delta_j] < [\Gamma:\Delta_{j+1}] < [\Gamma:\Delta_j]^{1+\delta}. \]
\end{thm}

\nid It follows from \cite{KN} that in general finitely generated groups fail to have an affirmative answer to much weaker versions of Question 2. Indeed, for any increasing function $f$, there exists a finitely generated group
$\Gamma$ such that for any family of relatively prime subgroups $\set{\Delta_j}$ there exist $j_0$ such that
\[ f([\Gamma:\Delta_{j_0}]) < [\Gamma:\Delta_{j_0+1}]. \]

\nid These results, aside from drawing analogies with the distribution of primes, are connected to a variant of subgroup growth that measures both relative index growth and the interplay between finite index subgroups. For this discussion we require some additional notation and terminology. Given a finitely generated, residually finite group $\Gamma$, let $\FI(\Gamma)=\set{\Delta_j}$ denote the set of finite index subgroups of $\Gamma$ enumerated by index. Similarly, $\NFI(\Gamma)$ denotes the subcollection of normal, finite index subgroups. For a set $S$ of $\Gamma$, we denote $S\smallsetminus\set{1}$ by $S^\bullet$. We define the \emph{divisibility function}
\[ \D_\Gamma\co \Gamma^\bullet \lra \N \]
by
\[ \D_\Gamma(\gamma) = \min \set{[\Gamma:\Delta_j] ~:~ \gamma \notin \Delta_j,~\Delta_j \in \FI(\Gamma)}. \]
We define the associated \emph{normal divisibility function} for normal, finite index subgroups in an identical way and denote it by $\D_\Gamma^\lhd$.\smallskip\smallskip

\nid  In this article we study the $L^1$--norm of these and more general divisibility functions. Rivin \cite{Rivin} studied a similar norm on free groups, while the asymptotic behavior of $L^\iny$--norm on metric $n$--balls was the focus of the articles \cite{Bou}, \cite{BM}, \cite{Hadad}, and \cite{Rivin}. More explicitly, we define the \emph{residual average} $\Ave(\Gamma)$ to be
\[ \Ave(\Gamma) = \int_{\wh{\Gamma}} \D_{\wh{\Gamma}}d\mu, \]
where $\wh{\Gamma}$ is the profinite completion of $\Gamma$, the measure $\mu$ is the Haar probability measure on $\wh{\Gamma}$, and $\D_{\wh{\Gamma}}$ is a continuous extension of $\D_\Gamma$ to $\wh{\Gamma}$. Similarly, the \emph{normal residual average} is
\[ \Ave^\lhd(\Gamma) = \int_{\wh{\Gamma}} \D_{\wh{\Gamma}}^\lhd d\mu. \]

\nid We will relate the finiteness of the above integrals to Questions 1 and 2---see Proposition \ref{Main}. In particular, an application of Theorem \ref{BertrandForArithmeticLattices} yields our next theorem.

\begin{thm}\label{GeneralMainTheorem}
If $\Gamma$ is a finitely generated linear group over $\C$, then $\Ave(\Gamma)$ is finite.
\end{thm}

\nid Using Theorem \ref{PowerBertrand}, we derive the following theorem.

\begin{thm}\label{GeneralMainNormalTheorem}
If $\Gamma$ is a finitely generated linear group over $\C$, then $\Ave^\lhd(\Gamma)$ is finite.
\end{thm}

\nid The key ingredient in the proof of Theorem \ref{GeneralMainTheorem} is Theorem \ref{BertrandForArithmeticLattices}; likewise Theorem \ref{PowerBertrand} is the key ingredient in the proof of Theorem \ref{GeneralMainNormalTheorem}. As the story leading to these connections is somewhat involved, we postpone a discussion of this here, hoping we have intrigued the reader sufficiently.\smallskip\smallskip

\nid We mention one complementary result to Theorem \ref{GeneralMainNormalTheorem}. Namely, if $\Gamma$ is either the first Grigorchuk group or $\SL(n,\Z_p)$ for $n>1$, then $\Ave^\lhd(\Gamma)$ is infinite---see Theorem \ref{FirstGrig}. These examples show that both linearity and finite generation are necessary in Theorem \ref{GeneralMainNormalTheorem}.

\paragraph{Article layout.}

In Section \ref{preliminary}, we develop the framework for residual averages. In Section \ref{GapConnection}, we relate finiteness of residual averages to gaps between subgroups. In Section \ref{SubquadraticSection}, we prove Theorems \ref{BertrandForArithmeticLattices} and \ref{GeneralMainTheorem}. In Section \ref{ElementarySection}, we prove Theorem \ref{PowerBertrand}. In Section \ref{ElementarySectionII}, we prove Theorem \ref{GeneralMainNormalTheorem}. In Section \ref{examplesOfFailure}, we discuss the examples above. We discuss integrating over other densities like the asymptotic, annular, and spherical densities in Section \ref{DensitySection}. Finally, in Section \ref{LarsenSect} we briefly mention connections this article has with certain zeta functions studied by Larsen \cite{Larsen}.

\paragraph{Acknowledgements.}

Foremost, we are extremely grateful to Benson Farb for his inspiration, comments, and guidance. We would also like to thank Emmanuel Breuillard for pointing out to us \cite[Lemma 3.1]{BG}. In addition, we would like to thank Misha Belolipetsky, Emmanuel Breuillard, Jordan Ellenberg, Skip Garibaldi, Tsachik Gelander, Fr\'{e}d\'{e}ric Haglund, Martin Kassabov, Jim Kelliher, Larsen Louder, Riad Masri, Alan Reid, and Tomasz Zamojski for numerous invaluable conversations on the material of this article. Finally, the second author was partially supported by an NSF postdoctoral fellowship.

\section{Preliminaries}\label{preliminary}

\nid We begin with a section that constructs a rigorous framework for residual averages.

\paragraph{1. Residual systems and divisibility functions.}

Throughout, $\Gamma$ will be an infinite, finitely generated, residually finite group. A collection $\Cal{F} = \set{\Delta_j}$ of finite index subgroups is a \emph{residual system} if
\[ \ba_{j=1}^\iny \Delta_j = \set{1}. \]
In addition, a residual system is \emph{relative prime} if each distinct pair $\Delta_j,\Delta_k$ in $\Cal{F}$ satisfies $\Gamma=\Delta_j\Delta_k$. It will be convenient to work with residual systems comprised of normal subgroups; we call such systems \emph{normal residual systems}.\smallskip\smallskip

\nid Associated to a residual system $\Cal{F}$ is the \emph{$\Cal{F}$--divisibility function}
\[ \D_\Cal{F}(\gamma) = \min\set{[\Gamma:\Delta_j]~:~\gamma \notin \Delta_j,~\Delta_j \in \Cal{F}}. \]

\paragraph{2. Profinite completions.}\label{ProfiniteCompletionsSection}

Let $\Cal{T}_{\textrm{pro}}$ be the profinite topology on $\Gamma$ given by declaring the finite index subgroups of $\Gamma$ to be a neighborhood basis for the identity element and by declaring left multiplication to be a homeomorphism. This topology is also the weak topology on the set of all homomorphisms of $\Gamma$ to finite groups, where we topologize the finite groups with the discrete topology $\Cal{T}_{\textrm{discrete}}$. \smallskip\smallskip

\nid There are several equivalent views of the profinite completion $\wh{\Gamma}$ of $\Gamma$. The profinite completion $\wh{\Gamma}$ is defined to be the inverse limit of the inverse limit system comprised of all finite quotients of $\Gamma$, where the finite quotients are equipped with the discrete topology. Consequently, $\wh{\Gamma}$ is a compact Hausdorff, topological group. We can also define the profinite completion to be the Cauchy completion of $\Gamma$ with respect to either a uniform structure on $\Cal{T}_{\textrm{pro}}$ (see \cite{Howes}) or via a metrization of $\Gamma$ (either can be used to define equivalent notions of Cauchy for sequences). For the latter, the finite generation of $\Gamma$ is required. We refer the reader to \cite{Wilson} for the general theory of profinite groups and profinite completions.\smallskip\smallskip

\nid Set $(\wh{\Gamma},\vp)$ to be the profinite completion of $\Gamma$ with associated continuous homomorphism \[ \vp\co \Gamma \lra \wh{\Gamma}. \]
The image of $\Gamma$ is dense and in the event that $\Gamma$ is residually finite, $\vp$ is injective. By work of Haar \cite{Haar}, since $\wh{\Gamma}$ is a compact topological group, $\wh{\Gamma}$ is endowed with a bi-invariant Haar measure $\mu$. Upon normalizing, we assume throughout that $\mu$ is the unique probability measure. \smallskip\smallskip

\nid Associated to every finite index, normal subgroup $\Delta$ in $\NFI(\Gamma)$ is an associated compact, open normal subgroup $\ol{\Delta}$ in $\wh{\Gamma}$ defined by taking the closure of $\vp(\Delta)$ in $\wh{\Gamma}$. The subgroup $\ol{\Delta}$ yields an extension
\[ \wh{q}_\Delta\co \wh{\Gamma} \lra \Gamma/\Delta \]
of the canonical epimorphism
\[ q_\Delta\co \Gamma \lra \Gamma/\Delta \]
that satisfies $q_\Delta = \wh{q}_\Delta \circ \vp$ (see \cite[Proposition 1.4.2]{Wilson}).\smallskip\smallskip

\nid We require the following elementary lemma for computational purposes.

\begin{lemma}\label{lemmaB}
Let $\Delta_1, \Delta_2$, be finite index subgroups of a residually finite group $\Gamma$. Then
\[ \ol{\Delta_1 \cap \Delta_2} = \ol \Delta_1 \cap \ol \Delta_2. \]
\end{lemma}

\nid Lemma \ref{lemmaB} follows from the well known fact that there exists a bijection between finite index, normal subgroups of $\Gamma$ and compact, open, normal subgroups of $\wh{\Gamma}$.

\paragraph{3. Extending divisibility functions to profinite completions.}

Having laid the groundwork for residual averages, we next extend $\D_\Gamma^\lhd$ to $\wh{\Gamma}$. The first step is the following lemma.

\begin{lemma}\label{ContinuityLemma}
\[ \D_\Gamma^\lhd\co (\Gamma^\bullet,\Cal{T}_{\textrm{pro}}) \lra (\N,\Cal{T}_{\textrm{discrete}}) \]
is continuous.
\end{lemma}

\nid For the proof of the lemma, recall that $\NFI(\Gamma) = \set{\Delta_j}$ is the collection of normal, finite index subgroups ordered such that for all $j$, the inequality
\[ [\Gamma:\Delta_j] \leq [\Gamma:\Delta_{j+1}] \]
holds. For each $j$, define
\begin{equation}\label{SubgroupDefinitions}
\Lambda_j = \ba_{\ell=1}^j \Delta_\ell, \quad \Gamma_j = \Lambda_{j-1}\Delta_j,
\end{equation}
and
\[ \NFI_n(\Gamma) = \set{\Delta \in \NFI(\Gamma)~:~ [\Gamma:\Delta]=n}. \]
We now prove Lemmas \ref{ContinuityLemma}.

\begin{pf}
We must show that for any subset $S$ of $\N$, the pullback $(\D_\Gamma^\lhd)^{-1}(S)$ is open. As
\[ (\D_\Gamma^\lhd)^{-1}(S) = \bu_{s \in S} (\D_\Gamma^\lhd)^{-1}(s), \]
it suffices to show that $(\D_\Gamma^\lhd)^{-1}(s)$ is open for all $s$ in $\N$. To this end, note that
\[ (\D_\Gamma^\lhd)^{-1}(s) = \bu_{\Delta_j \in \NFI_s(\Gamma)} \pr{\Delta_j^c \cap \Lambda_{j-1}}, \]
where $\Delta_j^c$ denotes $\Gamma \smallsetminus \Delta_j$. As $\Delta_j$, $\Lambda_j$ are both open and closed in the profinite topology, $(\D_\Gamma^\lhd)^{-1}(s)$ is open.
\end{pf}

\nid Since any sequence $\set{\gamma_j}$ in $\Gamma^\bullet$ which converges to the identity in $\Cal{T}_{\textrm{pro}}$ has the property that $\set{\D_{\Gamma}^\lhd(\gamma_j)}$ converges to infinity, we continuously extend $\D_{\Gamma}^\lhd$ to $\Gamma$ by declaring $\D_\Gamma^\lhd(1)$ to be infinity. The universal mapping property for profinite completions affords us with a unique, continuous extension
\[ \wh{\D_\Gamma^\lhd}\co \wh{\Gamma} \lra \N \cup \set{\infty} \]
such that the diagram
\[ \xymatrix{ \wh{\Gamma} \ar[rrr]^{\wh{\D_\Gamma^\lhd}} & & & \N \cup \set{\infty} \\ \Gamma \ar[u]^{\vp} \ar[rrru]_{\D_\Gamma^\lhd} & & &} \]
commutes. This yields a continuous map
\[ \wh{\D_\Gamma^\lhd}\co \wh{\Gamma} \lra \N \cup \set{\infty}, \]
where the target is giving the topology induced by the 1--point compactification of $\R$. As a result, we know that $\wh{\D_\Gamma^\lhd}$ is a measurable function and the integral
\[ \int_{\wh{\Gamma}} \wh{\D_\Gamma^\lhd} d\mu \]
is well defined. We define the \emph{normal residual average} of $\Gamma$ to be the integral
\begin{equation} \label{defAveEquality}
\Ave^\lhd(\Gamma) = \int_{\wh{\Gamma}} \wh{\D_\Gamma^\lhd} d\mu.
\end{equation}

\nid \textbf{Remark.}
By uniqueness of the extension and the universal mapping property, we see that $\wh{\D^\lhd_\Gamma} = \D^\lhd_{\wh \Gamma}$, where
\[ \D^\lhd_{\wh{\Gamma}}(\gamma) = \min \set{[\wh{\Gamma}:\ol{\Delta_j}]~:~ \ol{\Delta_j} \in \NFI(\wh{\Gamma})}. \]
For the remainder of the article, we will denote the extension by $\D^\lhd_{\wh\Gamma}$. \smallskip\smallskip

\nid The next lemma provides a sum formula for $\Ave^\lhd(\Gamma)$.

\begin{lemma} \label{defAvgLemma}
For $\Gamma,\Delta_j,$ and $\Lambda_j$ defined by (\ref{SubgroupDefinitions}), we have
\[ \int_{\hat \Gamma} \D_{\wh{\Gamma}}^\lhd d\mu = \sum_{j=1}^\iny [\Gamma:\Delta_j]\pr{\frac{[\Lambda_{j-1}:\Lambda_j]-1}{[\Lambda_{j-1}:\Lambda_j]}}\pr{\frac{1}{[\Gamma:\Lambda_{j-1}]}}. \]
\end{lemma}

\begin{pf}
\nid By the remark above, the definition of the Lebesgue integral, and properties of $\wh{\Gamma}$, we see that
\begin{align*}
\int_{\wh{\Gamma}} \D_{\wh{\Gamma}}^\lhd d\mu &= \sum_{n=1}^\iny n\mu((\D^\lhd_{\wh{\Gamma}})^{-1}(n)) \\
&= \sum_{n=1}^\iny \sum_{\Delta_j \in \NFI_n(\Gamma)} [\wh{\Gamma}:\ol{\Delta_j}]\mu\pr{\ol{\Lambda_{j-1}} \smallsetminus \ol{\Lambda_{j}}}.
\end{align*}
An elementary calculation yields
\[ \mu\pr{\ol{\Lambda_{j-1}} \smallsetminus \ol{\Lambda_{j}}} = \pr{\frac{[\ol{\Lambda_{j-1}}:\ol{\Lambda_j}]-1}{[\ol{\Lambda_{j-1}}:\ol{\Lambda_j}]}}\pr{\frac{1}{[\wh{\Gamma}:\ol{\Lambda_{j-1}}]}}. \]
By Lemma \ref{lemmaB}, we get
\[ \pr{\frac{[\ol{\Lambda_{j-1}}:\ol{\Lambda_j}]-1}{[\ol{\Lambda_{j-1}}:\ol{\Lambda_j}]}}\pr{\frac{1}{[\wh{\Gamma}:\ol{\Lambda_{j-1}}]}}
=\pr{\frac{[\Lambda_{j-1}:\Lambda_j]-1}{[\Lambda_{j-1}:\Lambda_j]}}\pr{\frac{1}{[\Gamma:\Lambda_{j-1}]}}, \]
finishing the proof of the lemma.
\end{pf}

\nid Replacing $\D_\Gamma^\lhd$ with $\D_\Gamma$ in the above discussion, yields a continuous extension $\D_{\wh{\Gamma}}$ of $\D_\Gamma$. We define the \emph{residual average} of $\Gamma$ to be
\[ \Ave(\Gamma) = \int_{\wh{\Gamma}} \D_{\wh{\Gamma}}d\mu. \]
Lemma \ref{defAvgLemma} also holds for $\D_\Gamma$. Finally, since $\D_{\wh{\Gamma}} \leq \D_{\wh{\Gamma}}^\lhd$, note that
\[ \Ave(\Gamma) \leq \Ave^\lhd(\Gamma). \]

\paragraph{4. General residual systems.}

For a normal, residual system $\Cal{F}$ on $\Gamma$, it need not be the case that the associated system $\wh{\Cal{F}} = \set{\ol{\Delta}}_{\Delta \in \Cal{F}}$ is a residual system. In fact, if we take the completion $\Cl_\Cal{F}(\Gamma)$ with respect to the system $\Cal{F}$, by the universal mapping property for the profinite completion, we have a surjective homomorphism
\[ \wh{\vp_\Cal{F}}\co\wh{\Gamma} \lra \Cl_\Cal{F}(\Gamma), \]
and $\ker \wh{\vp_\Cal{F}}$ measures the failure of $\wh{\Cal{F}}$ to be a residual system. On the other hand, we would like to define the $\Cal{F}$--residual average of the $\Cal{F}$--divisibility function $\D_\Cal{F}$. To this end, we define
\[ \wh{\D_\Cal{F}} = \D_{\Cl_\Cal{F}(\Gamma)} \circ \wh{\vp_\Cal{F}}, \]
where $\D_{\Cl_\Cal{F}(\Gamma)}$ is the extension of $\D_\Cal{F}$ to $\Cl_\Cal{F}(\Gamma)$. The construction of $\D_{\Cl_\Cal{F}(\Gamma)}$ is done in precisely the same way as the extension of $\D_\Gamma^\lhd$ to $\wh{\Gamma}$ was constructed above. We assert that
\begin{equation}\label{EqualIntegrals}
\int_{\wh{\Gamma}} \wh{\D_\Cal{F}}d\mu = \int_{\Cl_\Cal{F}(\Gamma)} \D_{\Cl_\Cal{F}(\Gamma)} d\mu_\Cal{F},
\end{equation}
where $\mu_\Cal{F}$ is the associated probability measure on $\Cl_\Cal{F}(\Gamma)$. To prove this equality, note that the map $\wh{\vp_{\Cal{F}}}$ has the property that for any measurable set $A$ in $\Cl_\Cal{F}(\Gamma)$, we have the equality
\[ \mu((\wh{\vp_\Cal{F}})^{-1}(A)) = \mu_\Cal{F}(A). \]
Consequently, for any measurable function $g$ on $\Cl_\Cal{F}(\Gamma)$, we have the equality
\[ \int_{\wh{\Gamma}} g \circ \wh{\vp_\Cal{F}} d\mu = \int_{\Cl_\Cal{F}(\Gamma)} gd\mu_\Cal{F}. \]
Taking $g$ to be $\D_{\Cl_\Cal{F}(\Gamma)}$ yields (\ref{EqualIntegrals}). In addition, the induced map
\[ \wh{\vp_\Cal{F}}^\star\co L^1(\Cl_\Cal{F}(\Gamma),\mu_\Cal{F}) \lra L^1(\wh{\Gamma},\mu) \]
is an isometric embedding, where for $g \in L^1(\Cl_\Cal{F}(\Gamma),\mu_\Cal{F})$, we have
\[ \wh{\vp_\Cal{F}}^\star(g) = g \circ \wh{\vp_\Cal{F}}. \]

\nid We define the $\Cal{F}$--residual average to be
\[ \Ave_\Cal{F}(\Gamma) = \int_{\wh{\Gamma}} \wh{\D_\Cal{F}}d\mu. \]
As before, we have
\[ \Ave_\Cal{F}(\Gamma) = \sum_{\Delta_j \in \Cal{F}} [\Gamma:\Delta_j]\pr{\frac{[\Lambda_{j-1}:\Lambda_j]-1}{[\Lambda_{j-1}:\Lambda_j]}}\pr{\frac{1}{[\Gamma:\Lambda_{j-1}]}}, \]
where $\Lambda_j$ is the intersection over the first $j$ subgroups in $\Cal{F}$.

\begin{lemma}\label{DominateLemma}
If $\Cal{F}$ is any residual system on $\Gamma$ for which $\Ave_\Cal{F}(\Gamma)$ is finite, then $\Ave(\Gamma)$ is finite.
\end{lemma}

\begin{lemma} \label{InfiniteQuotientCorollary}
Let $\psi\co \Gamma \to \Lambda$ be a surjective homomorphism of finitely generated, residually finite groups. If $\Lambda$ is infinite and $\Ave(\Lambda)$ is finite, then $\Ave(\Gamma)$ is finite. In addition, if $\Ave^\lhd(\Lambda)$ is finite, then $\Ave^\lhd(\Gamma)$ is finite.
\end{lemma}

\nid As the validity of Lemma \ref{DominateLemma} is clear, we only prove Lemma \ref{InfiniteQuotientCorollary}.

\begin{pf}[Proof of Lemma \ref{InfiniteQuotientCorollary}]
To begin, by the universal mapping property for profinite completions, the homomorphism $\psi$ has a continuous extension
\[ \wh{\psi}\co \wh{\Gamma} \lra \wh{\Lambda}. \]
The map $\wh{\psi}$ induces an isometric embedding
\[ \wh{\psi}^\star\co L^1(\wh{\Lambda},\mu_{\wh{\Lambda}}) \lra L^1(\wh{\Gamma},\mu_{\wh{\Gamma}}). \]
In particular, by hypothesis $\wh{\psi}^\star(\D_{\wh{\Lambda}}),\wh{\psi}^\star(\D_{\wh{\Lambda}}^\lhd)$ are in $L^1(\wh{\Gamma},\mu_{\wh{\Gamma}})$. The proof is completed by noting the inequalities
\[ \D_{\wh{\Gamma}} \leq \wh{\psi}^\star(\D_{\wh{\Lambda}}) \text{ and } \D_{\wh{\Gamma}}^\lhd \leq \wh{\psi}^\star(\D_{\wh{\Lambda}}^\lhd). \]
\end{pf}

\section{Gaps between subgroups and residual averages}\label{GapConnection}

\nid In this section, we relate the finiteness of residual averages to gaps between subgroups. This connection will be done via elementary group theory and Lemma \ref{defAvgLemma}.

\paragraph{1. Index sum formula.}

We begin by deriving a formula for $\Ave_\Cal{F}(\Gamma)$ in terms of indices of subgroups when $\Cal{F} = \set{\Delta_j}$ is a normal residual system. For $\Delta_j$, we have (see (\ref{SubgroupDefinitions}))
\begin{equation}\label{LatticeDiagram}
\xymatrix{ & \Gamma \ar@{-}[d]^{r_j} & \\ & \Gamma_j \ar@{-}[ld]_{t_j} \ar@{-}[rd]^{s_j} & \\ \Lambda_{j-1} \ar@{-}[rd]_{s_j} & & \Delta_j \ar@{-}[ld]^{t_j} \\ & \Lambda_j & }
\end{equation}
and from this we obtain
\[ \mu(\ol{\Lambda_{j-1}} \smallsetminus \ol{\Lambda_j}) = \pr{\frac{[\Lambda_{j-1}:\Lambda_j]-1}{[\Lambda_{j-1}:\Lambda_j]}}\pr{\frac{1}{[\Gamma:\Lambda_{j-1}]}} = \frac{s_j-1}{r_js_jt_j}. \]
As $[\Gamma:\Delta_j]= r_js_j$, we obtain our next lemma.

\begin{lemma}\label{SumFormula}
\[ \Ave_\Cal{F}(\Gamma) = \sum_{j=1}^\iny \frac{(s_j-1)}{t_j}. \]
\end{lemma}

\nid We say that a residual system $\Cal{F}$ is \emph{nested} if $\Delta_{j+1} < \Delta_j$ for all $j$. For a nested, normal residual system, from (\ref{LatticeDiagram}) it follows that $t_j=1$ for all $j$. In particular, by Lemma \ref{SumFormula}, we have the following lemma.

\begin{lemma}\label{NestedLemma}
If $\Cal{F}$ is a nested, normal residual system on $\Gamma$, then $\Ave_\Cal{F}(\Gamma)$ is infinite.
\end{lemma}

\nid We will use this in Section \ref{examplesOfFailure} to show the first Grigorchuk group and $\SL(n,\Z_p)$ have infinite normal residual averages.

\paragraph{2. Another index sum formula.}

We next aim to recursively compute the coefficients $t_j$. This is achieved with the following.

\begin{lemma}\label{RecursiveFormula}
Let $\Cal{F}= \set{\Delta_j}$ be a normal residual system. For all $j$, we have
\[ t_{j+1} = \frac{\prod_{\ell=1}^j s_\ell}{r_{j+1}}. \]
\end{lemma}

\begin{pf}
For this, consider the diagram:
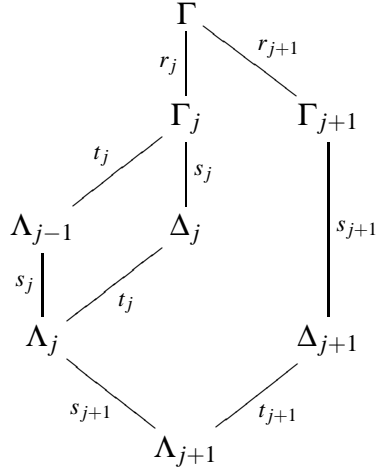
\[ \xymatrix{ & \Gamma \ar@{-}[d]_{r_j} \ar@{-}[rd]^{r_{j+1}} &  \\ & \Gamma_j \ar@{-}[ld]_{t_j} \ar@{-}[d]^{s_j} & \Gamma_{j+1} \ar@{-}[dd]^{s_{j+1}} \\ \Lambda_{j-1} \ar@{-}[d]_{s_j} & \Delta_j \ar@{-}[ld]^{t_j} & \\ \Lambda_j \ar@{-}[rd]_{s_{j+1}} & & \Delta_{j+1} \ar@{-}[ld]^{t_{j+1}} \\ & \Lambda_{j+1} & } \]
As we have two paths from $\Lambda_{j+1}$ to $\Gamma$, we see that
\[ [\Gamma:\Lambda_{j+1}] = t_{j+1}s_{j+1}r_{j+1} = s_{j+1}t_js_jr_j, \]
and thus
\begin{equation}\label{BigTree}
t_{j+1} = \frac{r_js_jt_j}{r_{j+1}}.
\end{equation}
To verify the formula for $t_{j+1}$, we utilize (\ref{BigTree}) via an inductive argument. For the base case, by definition, $r_1=1$ and $t_1=1$, and so
\[ t_2 = \frac{r_1s_1t_1}{r_2} = \frac{s_1}{r_2}. \]
We next assume now that the formula holds for $j$. From this assumption, we deduce the asserted formula:
\begin{align*}
t_{j+1} &= \frac{r_js_jt_j}{r_{j+1}} = \pr{\frac{r_js_j}{r_{j+1}}}\pr{\frac{\prod_{\ell=1}^{j-1} s_\ell}{r_j}} = \frac{\prod_{\ell=1}^j s_\ell}{r_{j+1}}.
\end{align*}
\end{pf}

\nid From Lemmas \ref{SumFormula} and \ref{RecursiveFormula}, we obtain another index sum formula for $\Ave_\Cal{F}(\Gamma)$:
\begin{equation}\label{SecondIntegralEquation}
\Ave_\Cal{F}(\Gamma) = \sum_{j=1}^\iny \frac{r_j(s_j-1)}{\prod_{\ell=1}^{j-1} s_\ell}.
\end{equation}

\paragraph{3. Index growth and convergence.} \label{indexGrowthAndConvergenceSection}

We are now in position to relate the finiteness of $\Ave(\Gamma)$ to gaps in subgroup growth. With $r_j,s_j$ given by (\ref{LatticeDiagram}), define the sequence
\[ \alpha_j(\Cal{F}) = \frac{r_{j+1}(s_{j+1}-1)}{r_js_j(s_j-1)}. \]
We say that $\Cal{F}$ as \emph{sub-quadratic index growth} if $\alpha_j(\Cal{F})<1$ for sufficiently large $j$. We say that $\Cal{F}$ as \emph{super-quadratic index growth} if $\alpha_j(\Cal{F})>1$ for sufficiently large $j$.

\begin{prop}\label{Main}
Let $\Gamma$ be a finitely generated, residually finite group and $\Cal{F}$ a normal residual system.
\begin{itemize}
\item[(a)]
If $\Cal{F}$ has sub-quadratic index growth, then $\Ave_\Cal{F}(\Gamma)$ is finite.
\item[(b)]
If $\Cal{F}$ has super-quadratic index growth, then $\Ave_\Cal{F}(\Gamma)$ is infinite.
\end{itemize}
\end{prop}

\begin{pf}
According to (\ref{SecondIntegralEquation}),
\[ \Ave_\Cal{F}(\Gamma) = \sum_{j=1}^\iny \frac{r_j(s_j-1)}{\prod_{\ell=1}^{j-1}s_\ell}, \]
where $r_j,s_j$ are given in (\ref{LatticeDiagram}). We see that the ratio of consecutive terms is given by
\begin{align*}
\frac{r_{j+1}(s_{j+1}-1)}{\prod_{\ell=1}^j s_\ell}\pr{\frac{r_j(s_j-1)}{\prod_{\ell=1}^{j-1} s_\ell}}^{-1} &= \frac{r_{j+1}(s_{j+1}-1)}{r_js_j(s_j-1)},
\end{align*}
and so both (a) and (b) follow from the Ratio Test.
\end{pf}

\nid We call a normal residual system $\Cal{F}$ \emph{prime} if $r_j=1$ for all $j$. One reason for this terminology is that prime systems have a property analogous to the Chinese Remainder Theorem:
\begin{equation}\label{ChineseRemainder}
\Gamma/ \Lambda_j = \bop_{\ell=1}^j \Gamma/\Delta_\ell.
\end{equation}
For such residual systems, the subgroups $\Delta_j$ have minimal overlap. Moreover, as $[\Gamma:\Delta_j]$ is unbounded, Proposition \ref{Main} reduces to studying the limit of the sequence
\[ \alpha_j(\Cal{F}) = \frac{[\Gamma:\Delta_{j+1}]}{[\Gamma:\Delta_j]^2}. \]
Therefore, the finiteness of $\Ave_\Cal{F}(\Gamma)$ depends only on the growth between consecutive indices for prime systems.\smallskip\smallskip

\nid For a normal residual system the coefficient $r_j$ measures the failure of (\ref{ChineseRemainder}) at the $j$th stage. Consequently, residual averages depend on the intersection of pairs and not just the growth of indices in general.\smallskip\smallskip

\nid For a non-normal residual system $\Cal{F}=\set{\Delta_j}$, we say that $\Cal{F}$ is \emph{prime} if for all $j$, we have
\[ [\Gamma:\Lambda_j] = \prod_{\ell=1}^j [\Gamma:\Delta_\ell]. \]
One case when this holds is when the indices $[\Gamma:\Delta_j]$ are pairwise relatively prime.

\paragraph{4. Residual averages on $\Z$.}

For $\Gamma =\Z$, divisibility functions have simple interpretations. For an integer $m$, we define three functions
\begin{align*}
\D(m) &= \min \set{n~:~ \GCD(m,n) < \min\set{m,n},~n>1} \\
\D_{\textrm{prime}}(m) &= \min \set{p~:~ \GCD(m,p)=1,~p \text{ prime}} \\
\D_p(m) &= \min\set{p^j~:~ \GCD(m,p^j)=1,~p\text{ a fixed prime}}
\end{align*}
We define the averages of these function as before and denote them $\Ave(\Z),\Ave_{\textrm{prime}}(\Z)$, and $\Ave_p(\Z)$, respectively. In the statement of our next result, $\set{p_j}$ denotes the set of primes ordered by cardinality.

\begin{prop}\label{MainIntegers}
We have the following formulas:
\begin{align*}
\Ave(\Z) &= \sum_{j=1}^\iny j\pr{1-\frac{\lcm(1,\dots,j-1)}{\lcm(1,\dots,j)}}\pr{\frac{1}{\lcm(1,\dots,j-1)}} \\
\Ave_{\textrm{prime}}(\Z) &= \sum_{j=1}^\iny \frac{p_j-1}{\prod_{\ell=1}^{j-1}p_\ell},\quad\quad  \Ave_p(\Z) = \sum_{j=1}^\iny (p-1)
\end{align*}
In particular, $\Ave(\Z),\Ave_{\textrm{prime}}(\Z)$ are finite while $\Ave_p(\Z)$ is not finite for any prime $p$.
\end{prop}

\nid This proposition follows easily using Bertrand's postulate, the Prime Number Theorem, and Proposition \ref{Main}. We also have (see Rivin \cite{Rivin} for the second series value):
\begin{align*}
\Ave(\Z) &\approx 2.787780456,\quad \Ave_{\textrm{prime}}(\Z) \approx 2.920050977.
\end{align*}

\section{Controlling gaps for linear groups}\label{SubquadraticSection}

In this section, we prove Theorem \ref{BertrandForArithmeticLattices}. Our proof splits into two cases depending on whether or not $\Gamma$ is virtually solvable.

\paragraph{1. Virtually solvable groups.}

When $\Gamma$ is virtually solvable, Theorem \ref{BertrandForArithmeticLattices} can be shown easily in a few different ways.

\begin{prop}\label{BertrandSolvableCase}
Let $\Gamma$ be a virtually solvable, finitely generated linear group over $\C$. Then there exists a constant $d$ and a family of finite index subgroups $\set{\Delta_j}$ of $\Gamma$ such that
\[ [\Gamma:\Delta_j] \leq [\Gamma:\Delta_{j+1}] \leq d[\Gamma:\Delta_j]. \]
In addition, there exists a finite index subgroup $\Gamma_0$ of $\Gamma$ such that $\set{\Delta_j}$ is a normal, prime family in $\Gamma_0$.
\end{prop}

\begin{pf}
It is well known that there exists a finite index subgroup $\Gamma_0$ of $\Gamma$ with the property that there exists a surjective homomorphism
\[ \vp\co \Gamma_0 \lra \Z. \]
Set $\Delta_j = \ker r_{p_j} \circ \vp$, where $\set{p_j}$ is the set of primes and
\[ r_{p_j}\co \Z \lra \B{F}_{p_j} \]
is reduction modulo $p_j$. This is a normal, prime family and by Bertrand's postulate, we have
\[ [\Gamma_0:\Delta_j] \leq [\Gamma_0:\Delta_{j+1}] \leq 2[\Gamma_0:\Delta_j]. \]
Viewed as subgroups of $\Gamma$, we obtain a family of finite index subgroups. Moreover,
\begin{align*}
[\Gamma_0:\Delta_j] &\leq [\Gamma_0:\Delta_{j+1}] \leq 2[\Gamma_0:\Delta_j] \\
[\Gamma:\Gamma_0][\Gamma_0:\Delta_j] &\leq [\Gamma:\Gamma_0][\Gamma_0:\Delta_{j+1}] \leq 2[\Gamma:\Gamma_0][\Gamma_0:\Delta_j] \\
[\Gamma:\Delta_j] &\leq [\Gamma:\Delta_{j+1}] \leq 2[\Gamma:\Delta_j],
\end{align*}
as needed.
\end{pf}

\nid An alternative to the above proof is to realize a finite index subgroup $\Gamma_0$ of $\Gamma$ as a finite index subgroup of $\B{S}(\Cal{O}_k)$ for a solvable, linear $k$--algebraic group $\B{S}$. We then use reduction homomorphism on $\B{S}(\Cal{O}_k)$ to produce a prime, residual system with the desired index gaps on $\Gamma_0$. Hence, in Proposition \ref{BertrandSolvableCase}, we also have
\[ \ba_{j=1}^\iny \Delta_j = 1. \]
Finally, we can also arrange it so that the subgroups $\Delta_j$ are normal in $\Gamma$ by appealing to the congruence subgroup property for $\B{S}(\Cal{O}_k)$ (see \cite{Chahal}). As neither of these properties are required in the sequel, we have opted to omit the details for these upgrades.

\paragraph{2. A simple example.}

For non-virtually solvable groups, we focus first on the case when $\Gamma = \B{G}(\Z)$ for a connected, simply connected, simple, linear $\Q$--algebraic group $\B{G}$. Before undertaking this endeavor, we present a simple, motivational example.\smallskip\smallskip

\nid \textbf{Example.} Let $\B{G}=\SL(n,\C)$ and $\Gamma = \SL(n,\Z)$. For each prime $p_j$, we have a surjective homomorphism
\[ r_{p_j}\co \SL(n,\Z) \lra \SL(n,\B{F}_{p_j}) \]
given by reducing coefficients modulo $p_j$. A simple computation shows that
\[ \abs{\SL(n,\B{F}_{p_j})} = \frac{\prod_{\ell=0}^{n-1} (p_j^n-p_j^\ell)}{p-1} = [\Gamma:\ker r_{p_j}]. \]
According to Bertrand's postulate, we know that
\[ \frac{[\Gamma:\ker r_{p_{j+1}}]}{[\Gamma:\ker r_{p_j}]} \leq \frac{(p_j-1)\prod_{\ell=0}^{n-1} ((2p_j)^n-(2p_j)^\ell)}{(2p_j-1)\prod_{\ell=0}^{n-1} (p_j^n-p_j^\ell)}. \]
We have a prime normal residual system $\{ \Delta_j \}$ by setting $\Delta_j = \ker r_{p_j}$.
Applying L'H\^{o}pital's rule, we obtain
\[ \lim_{j \to \iny} \frac{[\Gamma:\Delta_{j+1}]}{[\Gamma:\Delta_j]} \leq 2^{n^2-1}. \]
Therefore, for large values of $j$, we have
\[ [\Gamma:\ker r_{p_j}] \leq [\Gamma:\ker r_{p_{j+1}}] \leq 2^{n^2}[\Gamma:\ker r_{p_j}] = 2^{\dim \B{G} + 1}[\Gamma:\ker r_{p_j}]. \]

\paragraph{3. Integral points in simple linear algebraic groups.}

Using the method above, we prove the following proposition.

\begin{prop}\label{IntegralPointsCase}
Let $\B{G}$ be a connected, simply connected, simple, linear $\Q$--algebraic group and $\Gamma = \B{G}(\Z)$. Then there exists a constant $d$ and a family of normal, prime, finite index subgroups $\set{\Delta_j}$ of $\Gamma$ such that
\[ [\Gamma:\Delta_j] \leq [\Gamma:\Delta_{j+1}] \leq d[\Gamma:\Delta_j]. \]
In addition,
\[ \ba_{j=1}^\iny \Delta_j = 1. \]
\end{prop}

\nid In the proof, we write $f(p) \sim g(p)$ to mean $\lim\limits_{p \to \infty} \frac{f(p)}{g(p)} = 1$, where the limit is taken over $p$ in the intersection of the domains of $f$ and $g$.

\begin{pf}
According to the Strong Approximation Theorem (see \cite{Nori87}, \cite{Pink00}, or \cite{Weisfeiler84}), for all but finitely many primes $p$, the reduction modulo $p$ homomorphism maps $\Gamma$ surjectively onto $\B{G}(\B{F}_p)$. In addition, the kernels $\ker r_{p_j}$ yield a normal, prime, residual system on $\Gamma$. As $\abs{\B{G}(\B{F}_p)} \sim p^{\dim \B{G}}$ (see \cite[p. 123]{LS} or \cite[p. 131]{Steinberg}), by Bertrand's postulate, we see for consecutive primes $p_j,p_{j+1}$ that
\begin{align*}
\abs{\B{G}(\B{F}_{p_{j+1}})} &\sim p_{j+1}^{\dim \B{G}} \\
&\leq (2p_j)^{\dim \B{G}} \\
&= 2^{\dim \B{G}} p_j^{\dim \B{G}} \\
&\sim 2^{\dim \B{G}} \abs{\B{G}(\B{F}_{p_j})}.
\end{align*}
Thus, setting $d=2^{\dim \B{G}+1}$, we get the desired gap condition for the normal, prime family $\set{\ker r_{p_j}}$ for sufficiently large $j$.
\end{pf}

\paragraph{4. Non-solvable groups.}

\nid Using Proposition \ref{IntegralPointsCase} and the Lubotzky Alternative, we obtain the following corollary.

\begin{cor}\label{GeneralStrongBertrand}
Let $\Gamma$ be a finitely generated linear group over $\C$ that is not virtually solvable. Then there exists a constant $d$ and a family of finite index subgroups $\set{\Delta_j}$ of $\Gamma$ such that
\[ [\Gamma:\Delta_j] \leq [\Gamma:\Delta_{j+1}] \leq d[\Gamma:\Delta_j]. \]
In addition, there exists a finite index subgroup $\Gamma_0$ of $\Gamma$ such that $\set{\Delta_j}$ is a normal, prime family in $\Gamma_0$.
\end{cor}

\begin{pf}
As $\Gamma$ is not virtually solvable, the Lubotzky Alternative (see \cite[Theorem 16.4.12]{LS}) yields a finite index subgroup $\Gamma_0$ of $\Gamma$ and a representation
\[ \rho\co \Gamma_0 \lra \B{G}(\Z_{\Cal{P}_{\textrm{good}}}), \]
where $\B{G}$ is a connected, simply connected, simple, linear $\Q$--algebraic group, $\Cal{P}_{\textrm{good}}$ is a finite set of integral primes, and
\[ \Z_{\Cal{P}_{\textrm{good}}} = \bop_{p \notin \Cal{P}_{\textrm{good}}} \Z_p. \]
The content of the Lubotzky Alternative is that it guarantees that the Strong Approximation Theorem can be applied to $\rho(\Gamma_0)$. In particular, for all but finitely many primes (possibly more than $\Cal{P}_{\textrm{good}}$), reduction modulo $p$ maps $\Gamma_0$ onto $\B{G}(\B{F}_p)$. Taking $\Cal{F} = \set{\ker r_{p_j} \circ \rho}$, we obtain a normal prime family on $\Gamma_0$. Moreover, by Proposition \ref{IntegralPointsCase}, we have
\[ [\Gamma_0:\Delta_j] \leq [\Gamma_0:\Delta_{j+1}] \leq d[\Gamma_0:\Delta_j]. \]
Viewing these subgroups inside of $\Gamma$, as before, we see that
\[ [\Gamma:\Delta_j] \leq [\Gamma:\Delta_{j+1}] \leq d[\Gamma:\Delta_j] \]
still holds. Thus, we have a family of finite index subgroups on $\Gamma$ with the desired gap condition.
\end{pf}

\nid It could very well be that the subgroups $\Delta_j$ are not normal in $\Gamma$. Set
\[ \textrm{Core}(\Delta_j) = \ba_{\gamma \in \Gamma} \gamma^{-1}\Delta_j\gamma \]
to be the normal core of $\Delta_j$ in $\Gamma$. These subgroups yield the normal family $\set{\textrm{Core}(\Delta_j)}$ on $\Gamma$. We cannot ensure the gap condition for the indices since
\[ [\Gamma:\Delta_j] \leq [\Gamma:\textrm{Core}(\Delta_j)] \leq [\Gamma:\Delta_j]^{[\Gamma:\Gamma_0]}. \]
With regard to finiteness of normal residual averages, this is a problem.

\paragraph{5. The proof of Theorem \ref{BertrandForArithmeticLattices}.}

Let $\Gamma$ be a finitely generated linear group over $\C$. Recall that for Theorem \ref{BertrandForArithmeticLattices}, we must produce a family of finite index subgroups $\set{\Delta_j}$ such that
\[ [\Gamma:\Delta_j] \leq [\Gamma:\Delta_{j+1}] \leq d[\Gamma:\Delta_j] \]
for some constant $d$ and all $j$. If $\Gamma$ is virtually solvable, the existence is the content of Proposition \ref{BertrandSolvableCase}. Otherwise, the existence of such a family follows from Corollary \ref{GeneralStrongBertrand}. \qed

\paragraph{6. Proof of Theorem \ref{GeneralMainTheorem}.}

According to Theorem \ref{BertrandForArithmeticLattices}, there exists a normal, finite index subgroup $\Gamma_0$ of $\Gamma$ and a normal, prime family of finite index subgroups $\set{\Delta_j}$ on $\Gamma_0$ such that
\[ [\Gamma_0:\Delta_j] \leq [\Gamma_0:\Delta_{j+1}] \leq d[\Gamma_0:\Delta_j]. \]
Setting
\[ K = \ba_{j=1}^\iny \Delta_j, \]
and $\Lambda_0 = \Gamma_0/K$, the family $\set{\Delta_j}$ descends to a normal, prime residual system $\set{\Delta_j'}$ on $\Lambda_0$. Since
\[ \frac{[\Lambda_0:\Delta_{j+1}']}{[\Lambda_0:\Delta_j']^2} \leq \frac{d}{[\Lambda_0:\Delta_j']} \]
and $[\Lambda_0:\Delta_j']>d$ for all but finitely many $j$, by Proposition \ref{Main} and Lemma \ref{DominateLemma}, $\Ave(\Lambda_0)$ is finite. By construction, $K$ has infinite index, and thus by Lemma \ref{InfiniteQuotientCorollary}, $\Ave(\Gamma_0)$ is finite. We claim now that this implies that $\Ave(\Gamma)$ is finite. For this, we have
\begin{equation}\label{IntegralDecomposition}
\Ave(\Gamma) = \int_{\wh{\Gamma}} \D_{\wh{\Gamma}} d\mu = \int_{\wh{\Gamma}\smallsetminus \ol{\Gamma_0}} \D_{\wh{\Gamma}} d\mu + \int_{\ol{\Gamma_0}} \D_{\wh{\Gamma}}d\mu.
\end{equation}
The finiteness of integrals on the right hand side of (\ref{IntegralDecomposition}) can now be seen from the following two facts:
\[ \int_{\wh{\Gamma}\smallsetminus\ol{\Gamma_0}} \D_{\wh{\Gamma}} d\mu \leq [\Gamma:\Gamma_0]\mu(\wh{\Gamma} \smallsetminus \ol{\Gamma_0}) < \iny, \]
and
\[ \int_{\ol{\Gamma_0}} \D_{\wh{\Gamma}}d\mu \leq [\Gamma:\Gamma_0]\int_{\ol{\Gamma_0}} \D_{\ol{\Gamma_0}}d\mu = [\Gamma:\Gamma_0]\Ave(\Gamma_0) < \iny. \]
 \qed

\section{Proof of Theorem \ref{PowerBertrand}}\label{ElementarySection}

\nid In this section we prove Theorem \ref{PowerBertrand}.

\paragraph{1. The main proposition.}

\nid The following proposition is the main step in proving Theorem \ref{PowerBertrand}.

\begin{prop}\label{Main1}
Let $\Gamma$ be a finitely generated subgroup of $\GL(n,K)$ and $K/\Q$ a finite extension. Then for each $\delta>0$, there exists a normal residual system $\Cal{F}_\delta$ such that for all $j$,
\[ [\Gamma:\Delta_j] \leq [\Gamma:\Delta_{j+1}] \leq [\Gamma:\Delta_j]^{1+\delta}. \]
\end{prop}

\nid As the proof is somewhat involved, we summarize our strategy for the reader. Using the group $\SL(n,\Z)$ as a model, we take kernels of reduction homomorphisms
\[ r_\Fr{p}\co \Gamma \lra \GL(n,S/\Fr{p}) \]
for a particular subring $S$ in $K$ and prime ideals $\Fr{p}$ of $S$. Unlike the case of $\SL(n,\Z)$, we have no control here on the size of the index of $\ker r_\Fr{p}$. We circumvent this by instead taking reductions modulo $\Fr{p}_j^{k_j}$ for suitable powers $k_j$. The selection of these powers comprises the bulk of the proof. With regard to exposition, the difficulty is the interdependence of several quantities, each of which requires control for the selection of the powers $k_j$. With this in mind, in the proof below, we indicate the dependence of $x$ on $y$ by $x_y$. We hope this makes clear the dependence of each quantity on the others. These dependencies are important in both the proof of Proposition \ref{Main1} and of Theorem \ref{GeneralMainNormalTheorem}.

\begin{pf}[Proof of Proposition \ref{Main1}]
For a finitely generated subgroup $\Gamma$ of $\GL(n,K)$ and a finite generating set $\set{\gamma_m}$ of $\Gamma$, we define $S$ to be the ring generated by $\set{(\gamma_m)_{i,j}}$. By possibly enlarging $S$, we can assume that $\Cal{O}_K$ is contained in $S$, where $\Cal{O}_K$ is the ring of $K$--integers. As $\Gamma$ is finitely generated, the set of prime ideals of $S$ can be identified with a co-finite subset of $\Cal{P}_K$, where $\Cal{P}_K$ is the set of prime ideals of $\Cal{O}_K$. Specifically, we associate to $\Fr{p}$ in $\Cal{P}_K$ the ideal $S\Fr{p} = \Fr{p}_S$. For all but finitely many prime ideals $\Fr{p}$ in $\Cal{O}_K$, the ideal $\Fr{p}_S$ is a proper, prime ideal and
\[ S/\Fr{p}_S \cong \Cal{O}_K/\Fr{p} \cong \B{F}_{p^{s_{\Fr{p}}}}. \]
We denote the set of all prime ideals in $S$ by $\Cal{P}_S$. Excluding finitely many primes $p$, we select for each prime $p$ in $\Z$ a prime ideal $\Fr{p}_S$ in $S$ such that $S/\Fr{p}_S$ has characteristic $p$ and $\abs{S/\Fr{p}_S}$ is minimal. We denote this positive density subset of $\Cal{P}_S$ by $\Cal{P}_S^1$. According to the Cebotarev Density Theorem, we can pass to a positive density subset $\Cal{P}_S^2$ of $\Cal{P}_S^1$ such that for all $\Fr{p}_S$ in $\Cal{P}_S^2$, we have
\[ S/\Fr{p}_S = \B{F}_p. \]
We order $\Cal{P}_S^2 = \set{\Fr{p}_{S,j}}$ via the characteristic of the associated residue fields $\B{F}_p$.
By the Prime Number Theorem (see the proof of Lemma 2.4 in \cite{Bou}), the positive density of $\Cal{P}_S^2$ implies that there exists an integer $d$ such that
\begin{equation}\label{DensityInequality}
p_j \leq p_{j+1} \leq dp_j
\end{equation}
for all but finitely many $j$, where $p_j = |S/\Fr{p}_{S,j}|$.
Enlarging $d$ if necessary, we may assume that (\ref{DensityInequality}) holds for all $j$.
For each $\Fr{p}_{S,j}$ in $\Cal{P}_S^2$, we have the reduction homomorphism
\[ r_j\co \GL(n,S) \lra \GL(n,S/\Fr{p}_{S,j}) = \GL(n,\B{F}_{p_j}) \]
given by reducing the coefficients modulo the prime $\Fr{p}_{S,j}$. By our selection of $S$, the group $\Gamma$ is a subgroup of $\GL(n,S)$, and so the homomorphisms $r_j$ endow $\Gamma$ with a residual system $\Cal{F} = \set{\ker r_j \cap \Gamma}$. To see that this is a residual system, notice that for any $\gamma$ in $\Gamma^\bullet$, there are only finitely many prime ideals $\Fr{p}_{S,j}$ such that $\gamma$ resides in $\ker r_j$. In particular, we see that there are only finitely many $j$ such that $\Gamma$ is contained in $\ker r_j$, and we set $\Cal{P}_S^3$ to be the set of primes in $\Cal{P}_S^2$ for which $\Gamma$ is not a subgroup of $\ker r_j$. Set $\Cal{F}'=\set{\ker r_j \cap \Gamma}$ for $\Fr{p}_{S,j}$ in $\Cal{P}_S^3$. By construction, the image of $\Gamma$ under $r_j$ is a nontrivial subgroup of $\GL(n,\B{F}_{p_j})$ and thus has order $O_j$ where
\begin{equation}\label{OrderInequality}
1 < O_j \leq \prod_{\ell=0}^{n-1} (p_j^n - p_j^\ell) < p_j^{n^2}.
\end{equation}
To avoid controlling each $O_j$, we instead pass to deeper subgroups of $\Gamma \cap \ker r_j$ given by reduction modulo prime powers. If we reduce modulo the prime power $\Fr{p}_{S,j}^{k_j}$, we obtain the homomorphism
\[ r_{j,k_j}\co \GL(n,S) \lra \GL(n,S/\Fr{p}_{S,j}^{k_j}). \]
The image of $\Gamma$ under $r_{j,k_j}$ has order
\[ \abs{r_{j,k_j}(\Gamma)} = O_jp_j^{\ell_{j,k_j}} \]
This equality follows from the fact that for all $k>1$, we have (see \cite[Corollary 9.3]{Bass}, \cite[Ch. 9]{DDMS}, or the proof of Lemma 16.4.5 in \cite{LS})
\[ 1 \lra \Fr{gl}_n(\B{F}_{p_j}) \lra \GL(n,S/\Fr{p}_{S,j}^k) \lra \GL(n,S/\Fr{p}_{S,j}^{k-1}) \lra 1, \]
where $\Fr{gl}_n(\B{F}_{p_j})$ is the Lie algebra $\Fr{gl}_n$ of $\GL_n$ with coefficients in the finite field $\B{F}_{p_j}$. In particular, via induction, we have
\[ \abs{\GL(n,S/\Fr{p}_{S,j}^{k_j})} = p_j^{n^2(k_j-1)}\abs{\GL(n,\B{F}_{p_j})}. \]

\nid The associated residual system $\Cal{F}^\star = \set{\Gamma \cap \ker r_{j,k_j}}$ remains a normal residual system, and so we are reduced to finding a sequence of powers $\set{k_j}$ such that
\begin{equation}\label{RatioLimit}
\lim_{j \to \iny} \frac{[\Gamma:\Gamma \cap \ker r_{j+1,k_{j+1}}]}{[\Gamma:\Gamma \cap \ker r_{j,k_j}]^{1+\delta}} = \lim_{j \to \iny} \frac{O_{j+1}p_{j+1}^{\ell_{j+1,k_{j+1}}}}{O_j^{1+\delta}p_j^{(1+\delta)\ell_{j,k_j}}} < 1.
\end{equation}
We also require
\[ [\Gamma:\Delta_{j,k_j}] \leq [\Gamma:\Delta_{j+1,k_{j+1}}] \]
for large values of $j$. This yields our second desired inequality
\begin{equation}\label{BenIsStupid}
\frac{O_jp_j^{\ell_{j,k_j}}}{O_{j+1}p_{j+1}^{\ell_{j+1,k_{j+1}}}} < 1
\end{equation}
In order to achieve these inequalities, we select $N>(n^2)!$ and $C>4$. In addition, we may assume that $\delta<1/2$. We also select $0 < \eps < \delta$. In addition, let $j_{d,\eps}$ be such that for all $j>j_{d,\eps}$, we have $d<p_j^\eps$. \smallskip\smallskip

\nid For a fixed $j$, the sequence $\set{\ell_{j,k}}$ is non-decreasing and unbounded. That this sequence is non-decreasing followings from the short exact sequence. That this sequence is unbounded follows from
\[ \ba_{k=1}^\iny \ker r_{j,k} = 1. \]
Setting $\Delta_{j,k} = \ker r_{j,k}$, we seek a sequence $\set{k_j}$ such that
\[ [\Gamma: \Delta_{j,k_j}] < [\Gamma:\Delta_{j+1,k_{j+1}}] < [\Gamma:\Delta_{j,k}]^{1+\delta}. \]
To achieve these inequalities, we construct the sequence $k_j$ iteratively. To begin, select $k_1$ so that
$\ell_{1,k_1} > N + Cn^2$. Next, we select $k_2$ such that
\[ \ell_{1,k_1} + Cn^2 < \ell_{2,k_2} \leq \ell_{1,k_1} + (C+1)n^2. \]
As the sequence $\set{\ell_{2,k}}$ is non-decreasing and unbounded, there is a largest integer $i_2>1$ such that
\[ \ell_{2,i_2} \leq \ell_{1,k_1} + Cn^2. \]
In particular, we have
\[ \ell_{1,k_1} + Cn^2 < \ell_{2,i_2+1}. \]
Moreover, we have for all $j,k$ that
\[ \ell_{j,k+1} \leq \ell_{j,k} + n^2 \]
by the short exact sequence. In particular, we see that
\[ \ell_{1,k_1} + Cn^2 < \ell_{2,i_2+1} \leq \ell_{2,i_2} + n^2 \leq \ell_{1,k_1} + (C+1)n^2. \]
Setting $k_2 = i_2 + 1$, we also have
\[ \ell_{2,k_2} > \ell_{1,k_1} + Cn^2 > N + 2Cn^2. \]
Continuing iteratively, we produce a sequence $\set{k_j}$ such that for all $j$, we have
\[ \ell_{j,k_j} + Cn^2 < \ell_{j+1,k_{j+1}} \leq \ell_{j,k_j} + (C+1)n^2 \]
and
\[ \ell_{j,k_j} > N + Cjn^2. \]
We claim that this sequence achieves the above inequalities. We start with the inequality
\[ [\Gamma:\Delta_{j,k_j}] < [\Gamma:\Delta_{j+1,k_{j+1}}]. \]
We know that
\[ [\Gamma:\Delta_{j,k_j}] = O_jp_j^{\ell_{j,k_j}} < p_j^{n^2+\ell_{j,k_j}}. \]
We also have
\[ [\Gamma:\Delta_{j+1,k_{j+1}}] = O_{j+1}p_{j+1}^{\ell_{j+1,k_{j+1}}} > p_j^{\ell_{j+1,k_{j+1}}}. \]
Therefore, it suffices to have
\[ \ell_{j+1,k_{j+1}} > n^2 + \ell_{j,k_j}. \]
This inequality holds since $C>4$ and
\[ \ell_{j,k_j} + Cn^2 < \ell_{j+1,k_{j+1}}. \]
Next, we verify
\[ [\Gamma:\Delta_{j+1,k_{j+1}}] < [\Gamma:\Delta_{j,k}]^{1+\delta}. \]
Again, we have
\begin{align*}
[\Gamma:\Delta_{j+1,k_{j+1}}] &= O_{j+1}p_{j+1}^{\ell_{j+1,k_{j+1}}} < p_{j+1}^{n^2 + \ell_{j+1,k_{j+1}}} \\
&\leq d^{n^2+\ell_{j+1,k_{j+1}}}p_j^{n^2+\ell_{j+1,k_{j+1}}} \\
&< p_j^{\eps(n^2+\ell_{j+1,k_{j+1}})}p_j^{n^2+\ell_{j+1,k_{j+1}}} \\
&= p_j^{(\eps+1)(n^2+\ell_{j+1,k_{j+1}})} \\
&< p_j^{(\eps+1)(n^2 + \ell_{j,k_j} + (C+1)n^2)} \\
&= p_j^{(\eps+1)(\ell_{j,k_j} + (C+2)n^2)}.
\end{align*}
We also have
\[ [\Gamma:\Delta_{j,k_j}]^{1+\delta} = O_j^{1+\delta}p_j^{(1+\delta)\ell_{j,k_j}} \geq p_j^{(1+\delta)\ell_{j,k_j}}. \]
Thus, we require
\[ (1+\delta)\ell_{j,k_j} - (1+\eps)(\ell_{j,k_j} + (C+2)n^2) > 0. \]
By construction, we know that
\[ \ell_{j,k_j} > N + Cjn^2. \]
So we have
\begin{align*}
(1+\delta)\ell_{j,k_j} - (1+\eps)(\ell_{j,k_j} + (C+2)n^2) &= (\delta-\eps)\ell_{j,k_j} - (1+\eps)(C+2)n^2 \\
&> (\delta-\eps)[N + Cjn^2] - (1+\eps)(C+2)n^2.
\end{align*}
By selection, $\delta-\eps>0$ and $N,n,C,\eps$ are all constant. Thus, there exists $j_0 > j_{d,\eps}$ such that $p_j^\eps>d$ and
\[ (\delta-\eps)[N + Cjn^2] - (1+\eps)(C+2)n^2 > 1 \]
for all $j\geq j_0$. In total, we see that the normal residual system $\set{\Delta_{j,k_j}}_{j>j_0}$ satisfies the conditions needed for the proposition.
\end{pf}

\nid Without the Strong Approximation Theorem, we have very little control in the above proof. For instance, it is not clear that the sequence $\set{\ell_{j,k}}$ is strictly increasing for $k \geq M$, where $M$ is a constant that is independent of $j$. If the group $\Gamma$ is a cyclic subgroup $\innp{\alpha}$ of $\Cal{O}_K^\times$, the group of units in $\Cal{O}_K$, one would hope that for all but finitely many primes $\Fr{p}$, the order of $\alpha$ modulo $\Fr{p}^2$ is $a_\Fr{p}p$ where $q$ is the cardinality of the residue field $\Cal{O}_K/\Fr{p}$ and $a_\Fr{p}$ divides $q-1$. To put this problem into perspective, a Wieferich prime is a prime $p$ such that $p^2$ divides $2^{p-1} -1$. It has been conjectured that only finitely many Wieferich primes exist \cite{Wief} and also that infinitely many Wieferich primes exist \cite{Miri}. The above hope is analogous to the finiteness of Wieferich primes. Generalization of this problem are related to the ABC conjecture (see \cite{Silverman}).\smallskip\smallskip

\nid For a fixed prime ideal $\Fr{p}$, after the smallest power $k_\Fr{p}$ such that $\Fr{p}^{k_{\Fr{p}}}$ is not in primary decomposition of the ideal $(\alpha^{q-1} - 1)$, we do get strict growth on the prime powers $p$ occurring in the multiplicative order of $\alpha$ modulo $\Fr{p}^k$. This lack of control of $k_\Fr{p}$ makes estimates of the indices $s_j$ quite difficult. To prove Theorem \ref{GeneralMainNormalTheorem} using Theorem \ref{PowerBertrand}, the indices $s_j$ are precisely what requires control.

\paragraph{2. Proof of Theorem \ref{PowerBertrand}.}

We now prove Theorem \ref{PowerBertrand}. We start with a well known proposition.

\begin{prop}\label{RepresentationProposition}
If $\Gamma$ is a finitely generated, infinite linear group, then there exists an infinite representation of $\Gamma$ into $\GL(n,K)$ for some $n$ and finite extension $K/\Q$.
\end{prop}

\begin{pf}
If $\Gamma$ is virtually solvable, the conclusion of Proposition \ref{RepresentationProposition} follows from \cite[p.137]{Borel}. Otherwise, the conclusion follows, for instance, from \cite[Proposition 16.4.13]{LS} or \cite[Lemma 3.1]{BG}.
\end{pf}

\nid With Propositions \ref{Main1} and \ref{RepresentationProposition}, we now quickly derive Theorem \ref{PowerBertrand}.

\begin{pf}[Proof of Theorem \ref{PowerBertrand}]
According to Proposition \ref{RepresentationProposition}, there exists an infinite linear representation
\[ \rho\co \Gamma \lra \GL(n,K) \]
for some $n$ and $K/\Q$ finite. According to Proposition \ref{Main1}, for any $\delta>0$, there exists a normal residual system $\Cal{F}_\delta^\rho=\set{\Delta_j'}$ for $\rho(\Gamma)$ such that for each $\Delta_j',\Delta_{j+1}'$
\[ [\rho(\Gamma):\Delta_j'] \leq [\rho(\Gamma):\Delta_{j+1}'] \leq [\rho(\Gamma):\Delta_j']^{1+\delta}. \]
Setting
\[ \Cal{F}_\delta = \set{\Delta_j~:~ \Delta_j = \rho^{-1}(\Delta_j'), ~\Delta_j \in \Cal{F}_\delta^\rho}, \]
we see that
\[ [\Gamma:\Delta_j] \leq [\Gamma:\Delta_{j+1}] \leq [\Gamma:\Delta_j]^{1+\delta}. \]
That $\Cal{F}_\delta$ is a normal family of finite index subgroups of $\Gamma$ follows from elementary group theory.
\end{pf}

\section{Proof of Theorem \ref{GeneralMainNormalTheorem}}\label{ElementarySectionII}

\nid We are now ready to prove Theorem \ref{GeneralMainNormalTheorem}. Using the proof of Proposition \ref{Main1}, the main technical point is ensuring the coefficients $s_j = [\Lambda_{j-1}:\Lambda_j]$ are sufficiently large. We control these values via trivial estimates. The flexibility of the proof of Proposition \ref{Main1} allows us to use this growth condition to prove Theorem \ref{GeneralMainNormalTheorem} by appealing directly to Proposition \ref{Main}.

\begin{pf}[Proof of Theorem \ref{GeneralMainNormalTheorem}]
According to Proposition \ref{RepresentationProposition}, there exists an infinite representation
\[ \rho\co \Gamma\lra \GL(n,K) \]
for some $n$ and $K/\Q$ finite. By Lemma \ref{InfiniteQuotientCorollary}, it suffices to show that $\Ave^\lhd(\rho(\Gamma))$ is finite. For notational simplicity, set $\Lambda=\rho(\Gamma)$. Finally, set $S$ to be the coefficient ring of $\Lambda$.\smallskip\smallskip

\nid For each $\delta$, from the proof of Proposition \ref{Main1}, there exists a normal residual system $\Cal{F}_\delta$ on $\Lambda$ given by $\Delta_{j,k_j} = \Lambda \cap \ker r_{j,k_j}$, where
\[ r_{j,k_j}\co \GL(n,S) \lra \GL(n,S/\Fr{p}_{S,j}^{k_j}). \]
In addition, we have
\[ \abs{r_{j,k_j}(\Lambda)} = O_jp_j^{\ell_{j,k_j}} \]
where
\[ 1 \leq O_j < p_j^{n^2}. \]
We also have for constants $N>(n^2)!$ and $C>4$ that
\[ \ell_{j,k_j} > N+Cjn^2 \]
and
\[ \ell_{j,k_j} + Cn^2 < \ell_{j+1,k_{j+1}} \leq \ell_{j,k_j} + (C+1)n^2. \]
Finally, we will assume that $\delta<0$.\smallskip\smallskip

\nid To prove Theorem \ref{GeneralMainNormalTheorem}, we first give some trivial estimates for the size of $s_{j,k_j}$. for each $i < j$, we claim that the largest power of $p_j$ that divides $[\Lambda:\Delta_{i,k_i}]$ is $p_j^{n^2}$. To see this claim, note that if $p_j^m$ divides $O_ip_i^{\ell_{i,k_i}}$, since $p_i,p_j$ are distinct primes, $p_j^m$ must divide $O_i$. However, $O_i < p_i^{n^2}$ and $p_i < p_j$. Thus the claim follows. This claim thus shows that
\[ [\Delta_{i,k_i}: \Delta_{i,k_i} \cap \Delta_{j,k_j}] \geq p_j^{\ell_{j,k_j} - n^2}. \]
Taking this fact over all $i<j$, we see that
\[ [\Lambda_{j-1}: \Lambda_{j-1} \cap \Delta_{j,k_j}] \geq p_j^{\ell_{j,k_j} - (j-1)n^2}. \]
As the left hand side is $s_j$, we see that
\[ s_{j,k_j} \geq p_j^{\ell_{j,k_j}-(j-1)n^2}. \]

\nid We are now ready to show that $\Ave^\lhd(\Lambda)$ is finite. For this finiteness, by Proposition \ref{Main}, we must show directly
\[ \lim_{j \to \iny} \alpha_{j,k_j}(\Cal{F}_\delta) = \lim_{j \to \iny} \frac{r_{j+1,k_{j+1}}s_{j+1,k_{j+1}}-r_{j+1,k_{j+1}}}{r_{j,k_j}s_{j,k_j}(s_{j,k_j}-1)} = 0. \]
Thus, using Theorem \ref{PowerBertrand} and the definitions of $r_{j,k_j},s_{j,k_j}$ (see \ref{LatticeDiagram}), we see that
\begin{align*}
\frac{r_{j+1,k_{j+1}}s_{j+1,k_{j+1}} - r_{j+1,k_{j+1}}}{r_{j,k_j}s_{j,k_j}(s_{j,k_j}-1)} &< \frac{[\Lambda:\Delta_{j+1,k_{j+1}}]}{[\Lambda:\Delta_{j,k_j}](s_{j,k_j}-1)} \\
&\leq \frac{[\Lambda:\Delta_{j,k_j}]^\delta}{s_{j,k_j}-1}.
\end{align*}
As the term $s_{j,k_j}$ is increasing, we are reduced to showing
\[ \lim_{j\to \iny} \frac{[\Lambda:\Delta_{j,k_j}]^\delta}{s_{j,k_j}} = 0. \]
This limit is dealt with as follows.
\begin{align*}
\lim_{j\to \iny} \frac{[\Lambda:\Delta_{j,k_j}]^\delta}{s_{j,k_j}} &= \lim_{j\to \iny} \frac{O_j^\delta p_j^{\delta\ell_{j,k_j}}}{s_{j,k_j}} \\
&< \lim_{j\to \iny} \frac{p_j^{\delta n^2 + \delta \ell_{j,k_j}}}{p_j^{\ell_{j,k_j} - (j-1)n^2}}.
\end{align*}
Therefore, it suffices to have
\[ \ell_{j,k_j} - (j-1)n^2 - \delta n^2 - \delta\ell_{j,k_j} > 1. \]
To that end, we have
\begin{align*}
\ell_{j,k_j} - (j-1)n^2 - \delta n^2 - \delta\ell_{j,k_j} &= (1-\delta)\ell_{j,k_j} - [j+\delta-1]n^2 \\
&> \frac{1}{2}\ell_{j,k_j} - jn^2 \\
&> \frac{1}{2}[N + Cjn^2] - jn^2 \\
&= \frac{1}{2}N + \pr{\frac{Cj}{2}-j}n^2 \\
&> \frac{1}{2}N + jn^2 > n > 1.
\end{align*}
Hence, by Proposition \ref{Main}, $\Ave_{\Cal{F}_\delta}(\Lambda)$ is finite. Thus, by Lemma \ref{DominateLemma}, $\Ave^\lhd(\Lambda)$ is finite and so by Lemma \ref{InfiniteQuotientCorollary}, $\Ave^\lhd(\Gamma)$ is finite.
\end{pf}

\nid We obtain a different proof of Theorem \ref{GeneralMainTheorem} as $\Ave(\Gamma) \leq \Ave^\lhd(\Gamma)$.

\section{A pair of examples} \label{examplesOfFailure}

\nid In this short section, we show that neither the finite generation nor the linearity in Theorem \ref{GeneralMainTheorem} can be dropped. Recall that the first Grigorchuk group $\Gamma$ is a subgroup of automorphisms of a rooted binary tree (see \cite[Ch. VIII]{DeLa}) and is known to be finitely generated and residual finite but not linear. The groups $\SL(n,\Z_p)$, are linear and residually finite but not finitely generated for all $p$ and $n>1$.

\begin{thm}\label{FirstGrig}
The first Grigorchuk group $\Gamma$ has infinite $\Ave^\lhd$. The group $\SL(n, \Z_p)$ has infinite $\Ave^\lhd$ for all $p$ and $n>1$.
\end{thm}

\begin{pf}
For a normal, finite index subgroup $\Delta$ of $\Gamma$, the level of $\Delta$ is the largest $n$ such that $\Delta$ acts trivially on the $n$th level rooted binary tree. It is known that if $\Delta$ is level $n$, then $\Delta$ contains the kernel of the action on the $(n+3)$--level rooted binary tree (see \cite{CST}). This containment allows us to bound $\D_\Gamma^\lhd$ from below by a function $\D_1$ which has infinite $L^1$--norm. Specifically, set $\Delta_j$ to be the kernel of the induced map on the level $j$ rooted binary tree. We define
\[ \D_1(\gamma) = \begin{cases} 1, & \gamma \in \Delta_3^c \\ [\Gamma:\Delta_{j-3}], & \gamma \in \Delta_j^c \smallsetminus \Delta_{j-1}^c,~j>3. \end{cases} \]
Note that the function $\D_1$ is simply underestimating the value of $\D_\Gamma^\lhd$. In particular, we see that $\D_1 \leq \D_\Gamma^\lhd$. Moreover, we have by construction
\[ \int_{\wh{\Gamma}} \D_1d\mu = \pr{\frac{[\Gamma:\Delta_3]-1}{[\Gamma:\Delta_3]}} + \sum_{j=1}^\iny \frac{[\Gamma:\Delta_j]}{[\Gamma:\Delta_{j+2}]}\pr{\frac{[\Delta_{j+2}:\Delta_{j+3}]-1}{[\Delta_{j+2}:\Delta_{j+3}]}}. \]
Since the terms in latter series do not tend to zero, the above series diverges. As
\[ \int_{\wh{\Gamma}} \D_1d\mu \leq \Ave^\lhd(\Gamma), \]
we see that  $\Ave^\lhd(\Gamma)$ is infinite.\smallskip\smallskip

\nid For $\SL(n,\Z_p)$, the normal, finite index subgroups are all of the form $\ker r_j$ where $r_j$ is reduction modulo the $j$th power $\pi^j\Z_p$ of the uniformizer ideal $\pi \Z_p$. In particular, these subgroups are nested and hence by Lemma \ref{NestedLemma}, $\Ave^\lhd(\SL(n,\Z_p))$ is infinite.
\end{pf}

\nid Using the linear representations of $\Out(F_n)$ and $\Mod(S_g)$ obtained from acting on the first cohomology groups $H^1(F_n,\B{Z})$ and $H^1(\pi_1(S_g),\Z)$ respectively, we obtain from Lemma \ref{InfiniteQuotientCorollary} and Proposition \ref{Main1}, that $\Ave(\Out(F_n))$ and $\Ave(\Mod(S_g))$ are finite. Since for $n\geq 4$, $\Out(F_n)$ is not linear (see \cite{FP}), this provides an example of a non-linear, finitely generated, residually finite group with finite average. This shows that the finiteness of $\Ave(\Gamma)$ is not equivalent to the linearity of $\Gamma$. Finite generation is also not necessary, as $\SL(n,\wh{\Z})$ has finite average but is not finitely generated. Indeed, for any residually finite, finitely generated group $\Gamma$ with $\Ave(\Gamma)$ finite, we see that $\Ave(\wh{\Gamma})$ is finite.

\section{Averaging over other densities}\label{DensitySection}

\nid Throughout this article we study the average defined by integrating against the Haar measure in the profinite completion. In this section we average over densities, which are not measures in general. That is, for a finitely generated group $\Gamma$ and some density $\delta$, we define the \emph{average with respect to $\delta$} to be
\[ \Ave_\delta^\lhd (\Gamma) = \sum_{i=1}^\infty i \delta\left((\D^\lhd_\Gamma)^{-1}(i)\right). \]
For instance, given a group $\Gamma$ generated by a finite set $X$, the \emph{asymptotic density} of a subset $S$ of $\Gamma$ is defined to be
\[ \rho(S) = \limsup_{n \to \infty} \frac{| S \cap B_{X, \Gamma}(n) |}{|B_{\Gamma,X}(n) |}. \]
This density is never additive (see \cite{BV}, Example 2.2). There are other interesting notions of density, such as annular density, spherical density, and exponential density (see \cite{BV} and \cite{KRSS}). However when the asymptotic density exists as a limit, one can draw direct relationships between these other densities and the asymptotic density. In light of this fact, we focus on the asymptotic density.\smallskip\smallskip

\nid The asymptotic density is well-behaved when we restrict attention to $\Gamma$ satisfying:
\begin{equation} \label{BallGrowthEquation}
\lim_{n \to \infty} \frac{|B_{\Gamma,X}(n+1)|}{|B_{\Gamma,X}(n)|} = 1
\end{equation}
for some finite generating set $X$. In this case, by \cite[Section 2]{BV} the asymptotic density is both left and right invariant and we can replace the $\limsup$ with a limit in the definition of $\rho$ for finite index subgroups. In particular, for finite index subgroups, we have
\[ \rho(\Delta) = \frac{1}{[\Gamma:\Delta]}. \]
In tandem these aforementioned facts spawn our final result.

\begin{thm} \label{AsymptoticAverageTheorem}
Let $\Gamma$ be a residually finite group with finite generating set $X$ satisfying (\ref{BallGrowthEquation}). Then for all $n$ in $\N$,
\[ \mu\pr{(\D^\lhd_{\wh{\Gamma}})^{-1}(n)} = \rho((\D^\lhd_\Gamma(n))^{-1}). \]
In particular,
\[ \Ave^\lhd (\Gamma) = \Ave^\lhd_\rho(\Gamma). \]
\end{thm}

\begin{pf}
For the residual system $\NFI(\Gamma)$ with the subgroups $\Delta_j$ and $\Lambda_j$ given by (\ref{SubgroupDefinitions}), for any $i$ in $\N$, there exists $j,k$ in $\N$ such that $j \leq k$ and
\[ (\D^\lhd_\Gamma)^{-1}(i) = \Lambda_{j} \smallsetminus \Lambda_{k}. \]
By left invariance of $\rho$, we have
\[ \rho(\Lambda_{j} \smallsetminus \Lambda_k) = \rho(\Lambda_j) - \rho(\Lambda_k). \]
Hence, we have
\[ \rho(\Lambda_j \smallsetminus \Lambda_k) = \mu(\ol{\Lambda_j} \smallsetminus \ol{\Lambda_k}),\]
as desired.
\end{pf}

\section{Zeta functions associated to finite quotients}\label{LarsenSect}

\nid For an infinite group $\Gamma$, Larsen \cite{Larsen} studied the zeta function
\[ \zeta_\Gamma(s) = \sum_{i \in \Cal{I}(\Gamma)} \frac{1}{i^s}, \]
where
\[ \Cal{I}(\Gamma) = \set{i \in \N~:~ \NFI_n(\Gamma) \ne \es}. \]
He related the radius of convergence to $\min_\rho \dim\ol{\rho(\Gamma)}$, where
\[ \rho\co \Gamma \lra \GL(n,\C) \]
is a linear representation with infinite image and $\ol{\rho(\Gamma)}$ is the Zariski closure of $\rho(\Gamma)$. Visibly, convergence of $\zeta_\Gamma$ is related to gaps between successive indices and our work here provides some weak results on the radius of convergence. Specifically, we can give constants $s_\Gamma$ that ensure $\zeta_\Gamma(s)$ diverges provided $s \leq s_\Gamma$.



\noindent Department of Mathematics \\
University of Chicago \\
Chicago, IL 60637, USA \\
email: {\tt khalid@math.uchicago.edu}, {\tt dmcreyn@math.uchicago.edu}\\


\end{document}